\documentclass[11pt,letterpaper]{amsart}
\usepackage{amssymb, amscd, amsfonts, euler, epic, xy, charter}

\def\CC{{\mathbb C}}

\def\QQ{{\mathbb Q}} 
\def\RR{{\mathbb R}} 
\def\ZZ{{\mathbb Z}}

\def\half{\tfrac{1}{2}}
\def\G{{\Gamma}}

\def\dec{{\rm dec}}

\def\opp{{\rm op}}

\def\sep{{\rm sep}}

\def\ubold{\mathbf{u}}
\def\vbold{\mathbf{v}}
\def\Ybold{\mathbf{Y}}

\def\bs{\backslash}

\def\Acal{{\mathcal A}}

\def\Ccal{{\mathcal C}}
\def\Dcal{{\mathcal D}}
 
\def\Fcal{{\mathcal F}}

\def\Ical{{\mathcal I}}

\def\Ocal{{\mathcal O}}
\def\Scal{{\mathcal S}}
\def\Tcal{{\mathcal T}}
  
\def\Ucal{{\mathcal U}}

\def\la{\langle}
\def\ra{\rangle}

\newcommand\Hom{\operatorname{Hom}}

\newcommand\hgt{\operatorname{ht}}
\newcommand\HU{\operatorname{{\mathcal{HU}}}}

\newcommand\GL{\operatorname{GL}}

\newcommand\Sp{\operatorname{Sp}}

\newcommand\sign{\operatorname{sign}}

\newtheorem{theorem}{Theorem}[section]
\newtheorem{lemma}[theorem]{Lemma}
\newtheorem{proposition}[theorem]{Proposition}

\newtheorem{corollary}[theorem]{Corollary}

\theoremstyle{definition}

\theoremstyle{remark} 
\newtheorem{remark}[theorem]{Remark}

\title{Spherical complexes attached to symplectic lattices}
\author{Wilberd van der Kallen}
\author{Eduard Looijenga}
\email{W.vanderKallen@uu.nl, E.J.N.Looijenga@uu.nl}
\address{Mathematisch Instituut, Universiteit Utrecht, P.O.~Box 80.010, NL-3508 TA Utrecht (Nederland)}



\begin{document}\sloppy
\maketitle

\begin{abstract}To the integral symplectic group $\Sp(2g,\ZZ)$  we associate two posets of which we prove that they have 
the Cohen-Macaulay property. As an application we show that the locus of marked decomposable principally 
polarized abelian varieties in the Siegel space of genus $g$ has the homotopy type of a bouquet of 
$(g-2)$-spheres. This, in turn, implies that the rational homology of moduli space of (unmarked) principal 
polarized abelian varieties of genus $g$ modulo the decomposable ones vanishes in degree $\le g-2$.
Another application is an improved stability range for the homology of the symplectic groups over Euclidean 
rings. But the original motivation comes from envisaged applications to the homology of  groups of Torelli type. 

The proof  of our main result rests on a refined nerve theorem for posets that may have an interest in its own right.
\end{abstract}

\section*{Introduction}
\label{intro}
This paper is about  quasi-unimodular symplectic lattices, so let us begin with explaining that notion: 
a \emph{symplectic lattice} is a free abelian group $L$ of finite rank endowed with a symplectic form 
$(a,b)\in L\times L\mapsto \la a,b\ra$. It is said to be \emph{unimodular} if the associated map 
$a\in L\mapsto \la a,\;\ra\in \Hom (L,\ZZ)$ is an isomorphism; the rank of $L$ is then even and half 
that rank is called the \emph{genus} of $L$. Let us say that $L$ is \emph{quasi-unimodular} of genus 
$g$ if $L$ becomes unimodular of genus $g$ once we divide out by its radical.
For example, the intersection pairing on an oriented surface with finite first Betti number is 
quasi-unimodular with the genus equal to that of the surface. The unimodular sublattices of  $L$ make 
up a poset (with respect to inclusion) and our main result (Theorem \ref{thm:mc}) says that this poset 
is Cohen-Macaulay of dimension equal to the genus. 

However, for the applications we have in mind, the poset of unimodular decompositions of $L$ (where we 
now assume that $L$ is unimodular) is more relevant. We derive from our main result that this poset is 
also Cohen-Macaulay (Theorem \ref{thm:dec}). This has the following interesting 
consequence: consider in the Siegel space of genus $g$ the locus of decomposables, i.e., 
the locus that parametrizes marked principally polarized abelian varieties that are decomposable as 
polarized   varieties. This is a locally finite union of  symmetric submanifolds (each being isomorphic 
to a product of two Siegel spaces) and a nonempty  intersection of such submanifolds is also isomorphic 
to a product of  Siegel spaces (with perhaps more factors). 
A Siegel space is contractible  and so the locus of decomposables comes with a Leray covering by closed 
subsets. According to a classical result of Weil, the homotopy type of this locus is then 
that of the nerve of this covering. This nerve is just a  poset of unimodular decompositions and so we 
find (Corollary \ref{cor:dec}) that the  locus of decomposables has the homotopy type of a 
bouquet of $(g-2)$-spheres. (It would be desirable---and very interesting---to have a reasonable 
presentation of  the homology of this complex  in degree $g-2$ as a module over the symplectic group.) 
This implies that the rational homology of the pair $(\Acal_g,\Acal_{g,\dec})$,
where $\Acal_g$ is the moduli space of principally polarized abelian varieties and $\Acal_{g,\dec}$
parametrizes the decomposables, vanishes in degree $\le g-2$.

As an aside we observe that our results and arguments remain valid if we work over an arbitrary 
Euclidean ring  $R$ rather than $\ZZ$
 (so that $L$ is now a $R$-module). This makes it possible to improve Charney's  stability range for 
the homology groups 
of the groups $\{\Sp (2g,R)\}_{g\ge 0}$ by a unit or two (Theorem \ref{thm:improvedstab}). 
We may paraphrase this by saying that from her perspective, Euclidean rings behave as if they were fields.
 
Finally we show (Theorem \ref{thm:curvecx}) that for a closed orientable surface of genus $g$,  the 
\emph{separating curve complex} modulo the Torelli group has the same homotopy type as the locus of 
decomposables in Siegel space of genus $g$, so has  the homotopy type of a  bouquet of $(g-2)$-spheres. 
(One of us \cite{looij} recently proved that this also holds for separating curve complex itself.)
\\

We shall not review the individual sections, but we wish to point out the central role played by  
our Nerve Theorem \ref{thm:nerve}: while the first half of the nerve theorem may be familiar 
from Mirzaii-Van der Kallen (\cite[Thm.\ 4.3]{mk}); it is the second half that we believe is 
new and makes the theorem do the job that is needed here. We take the occasion to observe that 
its proof illustrates our belief that it is best when a proof of a statement  in any given category---here 
a homotopy category---does 
not leave that category. This means that we have expunged from the main argument all homology and the 
menagerie that usually accompanies it, such as the use of  spectral sequences, local systems and the 
Hurewicz theorem (and that remains so when we use it to derive the results that motivated this paper). 
The shortness of the proof of the Mirzaii-Van der Kallen part may be regarded as a testimony to the 
efficiency of this approach.

\section{The complex of unimodular sublattices}
\label{sec:1}Let $L$ be a free abelian group of finite rank endowed with a symplectic form 
$(a,b)\in L\times L\mapsto \la a,b\ra$. We say that $L$ is \emph{unimodular} if
the associated map $a\in L\mapsto \la a,\;\ra\in \Hom (L,\ZZ)$ is an isomorphism. Then $L$ has even rank  and half that rank is 
called the \emph{genus} of $g$.
We say $L$  is \emph{quasi-unimodular} if the associated map $a\in L\mapsto \bar a:=\la a,\;\ra\in \Hom (L,\ZZ)$ has a torsion 
free cokernel. So if $L_o\subset L$ denotes the kernel of this map (often referred to as the \emph{radical} of $L$), then the 
induced symplectic form on $\bar L:= L/L_o$ is unimodular; the \emph{genus} of $L$ is then by definition that of $\bar L$. 

We now fix a quasi-unimodular symplectic lattice $L$ of genus $g$.
The collection of unimodular  sublattices of $L$ (the trivial lattice included) make up a poset that we shall denote by $\Ucal (L)$. 
 Its standard height function is then given by the genus ($=$ half the rank).
A maximal chain is of the form ${0}=u_0\subsetneq u_1\subsetneq \cdots \subsetneq u_g$ and so $\dim \Ucal (L)=g$.

The main result of this article is 

\begin{theorem}\label{thm:mc}
The poset $\Ucal (L)$ is Cohen-Macaulay of dimension $g$.
\end{theorem}

We shall prove Theorem \ref{thm:mc} with induction on $g$. For $g=0$, $\Ucal (L)$ is a singleton and there is 
nothing to show. We assume the theorem verified for genera $<g$. 
Let us first isolate the two essential  cases needed for the induction step.

Since $\Ucal (L)$ has a minimal element (namely $0$), it is contractible  and hence  spherical of dimension $g$. 
If $u\in \Ucal (L)$, 
then $\Ucal (L)_{<u}=\Ucal (u)_{<u}$ is  spherical of dimension 
$g(u)-1$ for a similar reason. On the other hand,  $\Ucal (L)_{>u}$ can be identified with 
$\Ucal (u^\perp)_{>0}$ and hence if $g(u)>0$,
 this  is spherical of dimension $g-g(u)-1$ by our induction hypothesis. If $L$ is unimodular, then 
$L\in\Ucal (L)$ is a maximal element and
 then $\Ucal (L)_{>u}$ is even contractible. So the only case to deal with here is when $L$ is 
\emph{not} unimodular and $u=\{ 0\}$: 
we must show that $\Ucal (L)_{>0}$ is then $(g-2)$-connected. 

Similarly, if $u\subsetneq u'$ is an ordered pair in $\Ucal (L)$, then $u'':=u^\perp\cap u'$ has genus 
$g(u'')=g(u')-g(u)$ and so if $g(u')-g(u)<g$, then by our  induction hypothesis, 
$\Ucal (L)_{(u,u')}\cong\Ucal (u'')_{(0,u'')}$ is spherical of dimension $g(u'')-2$. 
So we may restrict to the case when $u'=0$ and 
 $g(u)=L$. In other words, we must show that when $L$ is unimodular, then  
$\Ucal (L)_{(0,L)}$ is $(g-3)$-connected. 
\\

Let $L$ be as above, i.e., quasi-unimodular of genus $g$. We denote the
poset  of isotropic sequences in $L$  that project to a partial basis of $\bar L$ by
$\Ical(L)$. For the induction step we need the following proposition, which is a special case of a result due to 
Barbara van den Berg (\cite[Prop.\ 1.6.1]{berg}). It is inspired by 
a similar result in the Utrecht  thesis (1979) of 
Maazen \cite{maazen}, which says that the poset of partial bases of a free finitely generated module over a 
Euclidean ring is 
Cohen-Macaulay. (See also the Appendix.) We here derive the result in question from Maazen's theorem.

\begin{proposition}\label{prop:st}
Let $L$ be a quasi-unimodular lattice of genus $g$. Then the poset $\Ical (L)$ of isotropic sequences that map to a 
partial basis of 
$\bar L$ is Cohen-Macaulay of dimension $g-1$.
\end{proposition}
\begin{proof} 
We first show that $\Ical (L)$  is spherical of dimension $g-1$.
Consider the poset $I(\bar L)$ of nonzero primitive isotropic sublattices of $\bar L$. According to a theorem of 
Solomon-Tits
 \cite[IV 5, Thm 2]{brown}, 
$I(\bar L)$ is Cohen-Macaulay of dimension $g-1$. So if $a\in I(\bar L)$ is of rank $k+1$, then $I(\bar L)_{>a}$ is 
spherical of 
dimension $g-2-k$. Let $f:\Ical (L)\to 
I(\bar L)$ be the poset map that assigns to an isotropic sequence in $L$ the span of its image in $\bar L$. 
Let $\tilde a$ be the preimage of $a$ under the projection $L\to \bar L$.
Then $f/a$ is the set 
of sequences  $(v_0,\dots ,v_s)$ in   $\tilde a$ that 
map to a partial basis of $a$. 
Or equivalently, if we first 
fix a basis $(w_1,\dots ,w_r)$ of $L_o$, that $(v_0,\dots ,v_s,w_1,\dots ,w_r)$ is a partial basis of $\tilde a$.  
According to  Maazen (\cite[Thm.\ III 4.2]{maazen} or Theorem \ref{thm:maazen} below),
this subposet is spherical of dimension $k$. %(Actually Maazen ignores fundamental groups,
%but his tools allow to treat fundamental groups also, 
%as illustrated by Ruth Charney in \cite{charney}.)
So condition  ($C^\opp$) of Corollary \ref{cor:conn} is fulfilled 
for $n=g-1$ and we conclude that $\Ical (L)$  is spherical of dimension $g-1$. 

We next verify the other properties needed for Cohen-Macaulayness.
Let $\vbold \in \Ical (L)$ have length $k+1$. Then $\Ical (L)_{<\vbold}$ is the poset of all proper subsequences of 
$\vbold$. This is essentially the boundary of a $k$-simplex and hence a $(k-1)$-sphere. 
For a similar reason, if $\vbold'$ is a subsequence of $\vbold$ of length $k'+1$, then 
$\Ical (L)_{>\vbold '}\cap \Ical (L)_{<\vbold}$ 
is spherical of dimension $k-k'-2$.

Finally, if $L'$ denotes the orthogonal complement of the span $I$ of $\vbold$ in $L$, then $L'$ is  
quasi-unimodular of genus $g-k-1$ with  radical $L_o\oplus I$ and $\Ical (L)_{>\vbold}$ consists of the sequences 
obtained by  shuffling
$\vbold$ with elements of $\Ical (L')$. By the above discussion combined with \cite[Cor.\ 1.7]{charney} for 
the shuffling, 
$\Ical (L)_{>\vbold}$ is therefore spherical of dimension  $g-k-2$.
\end{proof}

We use $\Ical (\bar L)$ to index full subcomplexes of  $\Ucal (L)_{>0}$ as follows. For a primitive vector $v\in \bar L$ 
denote by $L_v$ 
the set of $u\in L$ with $\la \bar u, v\ra=0$ and put $X_v:=\Ucal (L_v)$. 
For an isotropic sequence $\vbold =(v_0,\dots ,v_k)$, we 
put  $L_{\vbold}:=\cap_i L_{v_i}$  and $X_\vbold:=\cap_i X_{v_i}= X(L_{\vbold})$. The collection  
$\{\Ucal (L_\vbold)_{>0}\}_{\vbold\in \Ical (\bar L)}$ covers the poset of $u\in \Ucal (L)$ with 
$0< g(u)<g$ (which we shall denote 
by $\Ucal (L)_{(0,g)}$): if $u\in \Ucal (L)_{(0,g)}$, then there is a primitive isotropic $v\in\bar L$ perpendicular 
to $\bar u$  
and for any such $v$ we have $u\in \Ucal (L_v)_{>0}$. The fact that $\Ucal (L_\vbold)_{>0}=\emptyset$ if 
$\vbold$ has length $g$ will 
not bother us. 
Observe that if $\vbold'$ is a subsequence of $\vbold$, then $\Ucal (L_\vbold)_{>0}\subset \Ucal (L_{\vbold'})_{>0}$. 
\\
 
\begin{proof}[%
Proof 
of Theorem \ref{thm:mc} in case $L$ is unimodular]
We must show that $\Ucal (L)_{(0,g)}$ is $(g-3)$-connected. We do this by verifying the hypotheses (i) and (ii) of 
Theorem \ref{thm:nerve} for $X=\Ucal (L)_{(0,g)}$, for
$A:=\Ical (\bar L)$ indexing the collection of full subposets 
$\{X_\vbold:=\Ucal (L_\vbold)_{>0}\}_{\vbold\in \Ical (\bar L)}$ of 
$X$ and with $n=g-2$. We take for height functions the standard ones on these posets: $\hgt (u)=g(u)-1$ and 
$\hgt (\vbold)=|\vbold |-1$ (here $|\vbold|$ stands for the length of $\vbold$).

%Since $\Ical (\bar L_\vbold)$ is Cohen-Macaulay of dimension $g-1-\hgt (\vbold)$, ???
%$\Ical (\bar L_\vbold)$ is $(g-2-\hgt (\vbold))$-connected. So b
By our induction  hypothesis, the subset  
$X_\vbold=\Ucal (L_\vbold)_{>0}$ is $(g-2-\hgt (\vbold))$-connected. This verifies  \ref{thm:nerve}-(i).

Given  $u\in X$, then $X_{<u}=\Ucal(L)_{(0,u)}$ is $(g(u)-3)$-connected by induction. The poset 
$A_u$ of $\vbold\in\Ical (\bar L)_{>0}$  with $\bar u\perp \vbold$ is just
$\Ical (\bar u^\perp)$ and hence is $(g- g(u)-2)$-connected by Proposition \ref{prop:st}. 
Since $g- g(u)-2=(g-2)-\hgt (u)-1$,  \ref{thm:nerve}-(ii) is also satisfied.

Hence Theorem \ref{thm:nerve}  applies. Since $A$ is $(g-3)$-connected (it is even $(g-2)$-connected), it follows
that $X$ is  $(g-3)$-connected.  
\end{proof}

For the case when $L$ is not unimodular, we intend to invoke both halves of Theorem  \ref{thm:nerve}. We shall use the 
following elementary observation:

\begin{lemma}\label{lemma:elementary}
Let $\vbold=(v_0,\dots ,v_k)\in\Ical (\bar L)$.  If $u\in\Ucal (L)$ is of genus $k+1$ and such that 
$\vbold$ is contained in $\bar u$, then for every $u'\in  \Ucal (L_{\vbold})$, $u'+u$ is unimodular
of genus $g(u')+k+1$.
\end{lemma}
\begin{proof}
Choose a lift $e_i\in u$ of $v_i$ and extend $(e_0,\dots ,e_k)$ to a symplectic basis of $u$: 
$(e_0,\dots ,e_k;f_0,\dots ,f_k)$. For $u'$ as in the lemma, we have $e_i\perp u'$. Since $u'$ is unimodular,
there exists a $f'_i\in u'$ such that $f_i-f'_i$ is perpendicular to  $u'$. Then 
$(e_0,\dots ,e_k;f_0-f'_0,\dots ,f_k-f'_k)$
spans a unimodular lattice $\tilde u$ and $u+u'$ is the perpendicular direct sum of $\tilde u$ and $u'$.
So $u'+u$ is unimodular of genus $g(u')+k+1$.
\end{proof}

\begin{proof}[%
Proof 
of Theorem \ref{thm:mc} in case $L$ is not unimodular]
We wish to apply \ref{thm:nerve} to the case when $X=\Ucal (L)_{>0}$, but with the 
same collection of subposets $\{X_\vbold=\Ucal (L_\vbold)_{>0}\}_{\vbold\in \Ical (\bar L)}$ and value of $n$ 
(namely $g-2$) as 
in the unimodular case. The conditions (i) and (ii) are still satisfied for this value of $n$.
We also know that $A=\Ical (\bar L)$ is $(g-2)$-connected.

Let $\vbold=(v_0,\dots ,v_k)\in \Ical (\bar L)_k$.
For $i=0,\dots ,k$ we choose $u_i\in X$ of genus $1$ such that $v_i\in \bar u_i$ and 
$u_0,\dots ,u_k$ are pairwise perpendicular. For $\vbold'<\vbold$, we let $u_{\vbold'}\in  X$ be the span of the  
$u_{i}$ with $v_i\in\vbold-\vbold'$. Then $u_{\vbold'}\in  X$ is of genus $|\vbold|-|\vbold'|$ and is contained in 
$X_{\vbold'}$. According to 
Lemma \ref{lemma:elementary}, the sublattice  $u+u_{\vbold'}$ is unimodular
for any $u\in X_{\vbold-\vbold'}$. This is a fortiori so when
$u\in X_{\vbold}$, but since $u_{\vbold'}\in X_{\vbold'}$, we then in fact  have 
$u+u_{\vbold'}\in X_{\vbold'}$.

So condition \ref{thm:nerve}-(iii) is verified for  $s_{\vbold}(\vbold'):=u_{\vbold'}$ and 
$e_{\vbold}(\vbold', u):=u+u_{\vbold'}$. 
It remains to see that the resulting poset map
\[
\hat s_{\vbold} : (A_{<\vbold})^\opp\Join X_\vbold\to X
\]
is null-homotopic. The image $V_\vbold$ of $\hat s_{\vbold}$ is the union of $X_\vbold$ and the images of 
$s_{\vbold}$, $e_{\vbold}$.
Let $u_\emptyset$ denote the span of $u_0,\dots ,u_k$. Note that  $s_{\vbold}(\vbold')=u_{\vbold'}\le u_\emptyset$
and that if $u\in X_\vbold$, then $u+u_\emptyset\in X$ by Lemma \ref{lemma:elementary}.
So $\hat s_{\vbold}$  is homotopic to the constant map with value $u_\emptyset$ via the relations
$(A_{<\vbold})^\opp\Join X_\vbold\to V_\vbold\ni u\le u+u_\emptyset\ge u_\emptyset$.
\end{proof}

\section{Generalities on Posets}

We here collect some general results regarding posets that were used in the proofs of the previous section. 
A poset $X$ defines a simplicial complex with vertex set $X$ for which its $k$-simplices are chains 
$x_0<x_1<\cdots <x_k$. 
So the poset structure makes that every simplex has a natural order of its vertices, in particular,  has a natural 
orientation.
 The \emph{geometric realization} $|X|$ of $X$ is by definition that of the associated complex. So an element 
of $|X|$ is 
given by a 
function $\phi : X\to \RR_{\ge 0}$ whose support is a chain $x_0<x_1<\cdots <x_k$ and $\sum_i \phi (x_i)=1$.
We often allow ourselves the common  abuse of terminology when we say that $X$ enjoys a given (topological) property 
(such as 
connectedness or dimension), when we actually mean this to hold for $|X|$.

A $\ZZ$-valued function $f$ defined on the poset  $X$ is called a \emph{height function} if it is strictly increasing:  
$x<y$ implies $f(x)<f(y)$. If for every $x\in X$, $\dim (X_{\le x})$ (the supremum of the set of $n$ for which there 
exists a  chain $x_0<x_1<\cdots <x_n=x$ in $X$) is finite, then a standard choice for $f$ is $f(x):=\dim (X_{\le x})$. 
We call a height function \emph{bounded} if its image is contained in a finite interval. Given a height function 
$f$ on $X$,  we take $-f$ as height function on $X^\opp$.

%We use the convention that the reduced homology of a space $X$ is that of the mapping cylinder of the natural map from 
%$X$ to a 
%singleton shifted by one. So the reduced homology of the empty set is zero in all degrees $\not=-1$, 
%whereas $\tilde H_{-1}(\emptyset)=\ZZ$.
We use the convention that the dimension of the empty set is $-1$. 

Recall that a poset $X$ is said to be \emph{$n$-spherical} if its geometric realization is of dimension $n$ and 
$X$ is $(n-1)$-connected (we agree that the empty set is $(-1)$-spherical). It is said to be \emph{Cohen-Macaulay} of 
dimension $n$ if in addition
\begin{enumerate}
\item[(i)] for every $x\in X$, $X_{<x}$ resp.\  $X_{>x}$  is spherical of dimension  
$\dim X_{\le x}-1$ resp.\ $n-1-\dim X_{\le x}$ and 
\item[(ii)] for every ordered pair $x<y$ in $X$, $X_{>x}\cap X_{<y}$ is spherical of dimension
$\dim X_{\le y}-\dim X_{\le x}-2$.
\end{enumerate}
%We must assume some familiarity with \cite{quillen}.
If $X$,$Y$ are posets we define its join $X* Y$ as in \cite[1.8]{quillen} to be the disjoint union of 
$X$ and $Y$ equipped with the ordering which agrees with the given orderings on $X$ and $Y$ and which is such 
that any element of $X$ is less than any element of $Y$. Note that this is asymmetric in $X$ and $Y$.  
It is clear that with this definition an iterated join $X_1* X_2* \cdots * X_n$ has as underlying set the 
disjoint union of $X_1,\dots ,X_n$ with the given partial order on any piece $X_i$ and with any element of 
$X_i$ dominating $X_1\cup\cdots \cup X_{i-1}$.
In the case of two posets $X,Y$ we will also make use of the \emph{thick join}  $X\Join Y$ which is symmetric in 
$X$ and $Y$. It is
defined as follows. Letting $X\times Y$ denote  the product of $X$ and $Y$ in the category of posets, 
then $X\Join Y$ is the disjoint union of $X$, $Y$, $X\times Y$ equipped with the ordering which
is obtained by adding to the given orderings on $X$, $Y$, $X\times Y$ the relations $x<(x,y)>y$ for $x\in X$, $y\in Y$.
The poset map $h(X,Y):X\Join Y\to X* Y$ which maps $x$ to $x$, $(x,y)$ to $y$, $y$ to $y$, is a homotopy equivalence 
\cite[Prop.\ II 1.2]{maazen}.
%Indeed Quillen's Theorem A \cite[Prop.\ 1.6]{quillen} applies. 
If $X$ is $d$-connected and $Y$ is $d'$-connected, then $X* Y$ 
is $d+d'+2$-connected (and hence $X\Join Y$ is too).

Let $f:X\to Y$ be a map of posets. 
Recall that $f/y=\{x\in X\mid fx\leq y\}$ and $y\backslash f=\{x\in X\mid fx\geq y\}$.
We define the \emph{mapping cylinder} of $f$, $M(f)$,  to be the disjoint union of 
$X$ and $Y$ equipped with the ordering which
agrees with the given orderings on $X$ and $Y$ and which is such that for
$x\in X$, $y\in Y$, one has $x>y$ if and only if  $fx\geq y$. Observe that $Y$ is a deformation retract of $M(f)$.
Dually we define $M^\opp(f):=M(f^\opp)^\opp$, where $f^\opp:X^\opp\to Y^\opp$. (As a set map $f^\opp$ equals $f$.)
Again $Y$ is a deformation retract, but in $M^\opp(f)$ any element
of $x\in X$ is less than  $fx\in Y$.
If $Y$ is a %finite dimensional 
poset equipped with a %standard 
height function $\hgt$ we let
$M(f,\leq k)$ denote the disjoint union of $X$ and $Y_{\leq k}:=\{ y\in Y \mid \hgt(y)\leq k\}$
with the ordering induced from $M(f)$. 

Let us also define the \emph{mapping cone} of $f$, although this notion will play a less prominent 
role in this paper:  it is the disjoint union $C(f)$ of  a singleton $\{v_f\}$ and $M(f)$ whose ordering 
extends the one on $M(f)$ and for which $v_f$ is smaller than every point of $X$ and is incomparable 
with the points of $Y$.

%With notations as in \cite{quillen} we have in fact

\begin{proposition}\label{prop:grow}
Let $f:X\to Y$ be a map of  posets and let $\hgt$ be a bounded height function on $Y$.
Then $M(f,\leq k)=X$ for $k\ll 0$, $M(f,\leq k)= M(f)\sim Y$ for $k\gg0$. Furthermore, 
for every $y\in Y$ of height $k$,  the link of $y$ in $M(f,\leq k-1)$ can be identified with
$Y_{<y}*(y\backslash f)$, and so $M(f,\leq k)$ is obtained  from $M(f,\leq k-1)$ by putting for every 
$y$ with $\hgt(y)=k$ a cone over $Y_{<y}*(y\backslash f)$. If in addition,
$Y_{<y}* (y\backslash f)$ is $(n-1)$-connected  for all $y\in Y$, then   
$f$ is $n$-connected (which means that the pair  $(M(f), X)$ is $n$-connected).
%In particular, $X$ is $(n-1)$-connected if and only if $Y$ is $(n-1)$-connected. 
%Moreover $\tilde H_n(X)\to \tilde H_n(Y)$ is then onto.
\end{proposition}
\begin{proof}
The first assertion is clear and so is the second. If $Y_{<y}*(y\backslash f)$ is $(n-1)$-connected for all $y\in Y$, 
then $|M(f)|$ is gotten, up to homotopy, from $|X|$ by successive attachment of cells of dimension $\ge n+1$ and so the 
pair $(M(f),X)$ 
is $n$-connected.
\end{proof}
% Hatcher, Algebraic Topology, Chapter 0 page 11 implies that if one attaches something to a CW subspace $A$ of a CW space $X$
% so as to make it contractible, then this amounts to collapsing that subspace, up to homotopy equivalence. 
% And to kill the homotopy groups of an $(n-1)$-connected subspace $A$ you do not need to attach a full cone $CA$. 
% It can be done by successive attachment of cells of dimension $\ge n+1$. Say the result is a pair $(Z,A)$.
% Now to see that $(X\cup_A CA,X)$ is $n$-connected, one may look at $(X\cup_A Z,X)$ instead, as both problems are equivalent
% to $X\to X/A$ being $n$-connected.
% For an empty subspace $A$ one needs a separate argument.
% Alternative argument for homology: The mapping cone of space-->space+(cone on link) is homotopy equivalent to suspension of
% link.

We get the following slight variation on Theorem 9.1.\ of Quillen
\cite{quillen}. Compare also \cite[Thm 3.8]{mk}. 

\begin{corollary}\label{cor:conn}
Let $f:X\to Y$ be a map of  posets. Assume $Y$ is endowed with a bounded
height function $\hgt$ with the property that for some integer $n$ and some set map $t:Y\to \ZZ$ %(think $t=\hgt$) 
one of  the following is true:
\begin{enumerate}
\item[($C$)] for every $y\in Y$, $Y_{<y}$ is $(t(y)-2)$-connected and $y\bs f$ is $(n-t(y)-1)$-connected or dually,
\item[($C^\opp$)] for every $y\in Y$, $Y_{>y}$ is $(n-t(y)-2)$-connected and $f/y$ is $(t(y)-1)$-connected.
\end{enumerate}
Then $f$ is $n$-connected.
% and we have for every integer $k$ a
%natural exact sequence 
%\begin{align*}
%\bigoplus_{\hgt(y)=k}\tilde H_n(Y_{<y}*(y\backslash f))\to\tilde 
%H_n(M(f,\leq k-1))\to \tilde H_n(M(f,\leq k))\to 0,\\%\tag{$E_k$}\\
%\bigoplus_{\hgt(y)=k}\tilde H_n((f/y)*Y_{>y})\to\tilde 
%H_n(M^\opp(f,\geq k+1))\to \tilde H_n(M^\opp(f,\geq k))\to 0 %\tag{$E^\opp_k$}.
%\end{align*}
%in the  case ($C$) resp.\  ($C^\opp$). 
%Suppose this connectivity property holds.
%denote by $F_k\subset \tilde H_n(X)$ the image of 
%$\tilde H_n(f^{-1}Y_{\le k})\to \tilde H_n(X)$ so that we have an increasing filtration 
%\[
%0=F_{-1}\subset F_0\subset\cdots \subset F_{n}\subset F_{n+1}=\tilde H_n(X).
%\]
\end{corollary}
\begin{proof}
In case ($C$), first observe that $Y_{<y}*(y\backslash f)$ is $(t(y)-2)+(n-t(y)-1)+2$-connected,
hence $n-1$-connected, so that the result follows from Proposition \ref{prop:grow}. In case  ($C^\opp$), pass to the opposite.
\end{proof}

In many cases of interest, $t=\hgt$ does the job. The same is true for the following nerve theorem, which  is the technical 
result that makes the proof of Main Theorem \ref{thm:mc} possible. 
The first half of the nerve theorem is familiar from Mirzaii-Van der Kallen (\cite[Thm.\ 4.3]{mk}). 
We reprove it here to illustrate the efficiency of the present setup.
The second half is our key advance.
%We leave it to the reader to formulate a common generalization.

\begin{theorem}[Nerve Theorem for Posets]\label{thm:nerve} 
Let $X$ and $A$ be  posets both endowed with bounded height functions. Assume that $A^\opp$ labels \emph{full} 
subposets of $X$ in the sense  
for every $a\in A$ we are given a subposet $X_a$ with the property that  if $x\in X_a$, then also any $y<x$ is 
in $X_a$ and  $a<b$ 
implies $X_a\supset X_b$. Let $n$ be an integer such that for some set maps 
 $t_X:X\to \ZZ$  and  $t_A:A\to \ZZ$
\begin{enumerate}%\setlength{\itemindent}{.3em}
\item[(i)]  for every $a\in A$,  $A_{<a}$ is $(t_A(a)-2)$-connected, 
$X_a$ is $(n-t_A(a)-1)$-connected and 
\item[(ii)] for every $x\in X$,  $X_{<x}$ is $(t_X(x)-2)$-connected and the full subcomplex $A_x$ of $A$ spanned by the 
$a\in A$ 
with $x\in X_a$ is $(n-t_X(x)-1)$-connected.
\end{enumerate}
Then $A$ is  $(n-1)$-connected if and only if $X$ is. 

If in addition to (i) and (ii) there exist for every $a\in A$, poset maps 
$s_a: (A_{<a})^\opp\to X$ and $e_a: (A_{<a})^\opp\times X_a\to X$  such that 
\begin{enumerate}%\setlength{\itemindent}{.3em}
\item[(iii)] if $b<a$, then $s_a(b)\in X_b$ and for all $x\in X_a$ one has $e_a(b,x)\in X_b$ and 
$s_a(b)\le e_a(b,x)\ge x$, and
\item[(iv)] the resulting poset map $\hat s_a:(A_{<a})^\opp\Join X_a\to X$ is null-homotopic,
\end{enumerate}
then there even exists an $n$-connected map $|A|\to |X|$ (so that
$X$ is $n$-connected when $A$ is).
\end{theorem}
\begin{proof} Denote by 
$Z\subset A^\opp\times X$ the subset of pairs $(a,x)\in A\times X$ with $x\in X_a$ with the induced partial order. 
%by $Z_k\subset Z$ the subposet ofv $Z$ defined by
%$\hgt (a)\le k$. 
We have poset projections
\[
f: Z^\opp\to A\; \text{  and }\; g: Z \to  X. 
\]
We now divide the proof in a number of steps. Since $|Z^\opp|= |Z|$, the first assertion of the theorem follows from:
\\

\emph{Step 1: Both $f$ and $g$ are $n$-connected.} For every $a\in A$, we have that
$a\bs f=\{ (b,y)\in Z\, |\, b\ge a, y\in X_b\}$ contains $\{ a\}\times X_a$ as deformation retract with retraction given by 
$(b,y)\mapsto (a,y)$. So by (i), $a\bs f$ is $(n-\hgt (a)-1)$-connected. By that same assumption, $A_{<a}$ is  
$(\hgt (a)-2)$-connected. According to Corollary \ref{cor:conn}, $f$ is the $n$-connected.
A similar argument applied to $g: Z\to X$, with $g/x$ instead of $x\backslash g$, yields that $g$ is $n$-connected.  
\\

{}From now on we assume that the conditions (iii) and (iv) are satisfied. 
\\

\emph{Step 2: Construction of a diagram, commutative up to homotopy.}
Let $a\in A$ and put $k=\hgt(a)$.
We shall need the following diagram of poset maps
\[
\begin{CD}
((A_{<a})^\opp\Join X_a)^\opp @>{h^\opp}>{\sim}>  A_{<a} *(X_a)^\opp  @>{1*j^\opp}>{\sim}> A_{<a}*(a\bs f)^\opp\\
@V{\tilde s_a^\opp}VV  @. @VV{\text{attaching map}}V\\
Z^\opp @>>> M(f,\leq k-1) @= M(f,\leq k-1).
\end{CD}
\]
Here the top horizontal maps are  homotopy equivalences: the first map is the opposite of the natural map 
$h: X_a\Join (A_{<a})^\opp\to X_a *(A_{<a})^\opp$ encountered before (it sends $b\in A_{<a}$ to $b$, $x\in X_a$ to $x$ 
and $(b,x)$ to $b$) and the second map
is the join of the identity map of $A_{<a}$ and the opposite of the homotopy equivalence 
$j: x\in (X_a)\mapsto (a,x)\in (a\bs f)$. The first  lower horizontal map is the inclusion.
The vertical map on the left is the opposite of the poset map 
\begin{multline*}
\tilde s_a: (A_{<a})^\opp\Join X_a\to Z,\\
\tilde s_a \big(b<(b,x)>x\big)=\big((b,s_a(b))\le (b,e_a(b,x))\ge(a,x)\big),
\end{multline*}
where $b\in A_{<a}$ and $x\in X_a$. Notice that $g\tilde s_a=\hat s_a$.
The vertical map on the right is the embedding of the link of $a\in M(f,\le k)$ in  $M(f,\leq k-1)$ 
that appears in Proposition \ref{prop:grow} (it is given by $A_{<a}\subset A$ and $(a\bs f)^\opp\subset Z^\opp$).

We check that the diagram  is homotopy commutative.
The two maps from $((A_{<a})^\opp\Join X_a)^\opp$ to $M(f,\leq k-1)$ (whose underlying set is
$A_{\le k-1}\sqcup Z^\opp $) are given by
\[
\big(b>(b,x)<x \big)\mapsto 
\begin{cases}
\big((b,s_a(b))\ge (b,e_a(b,x))<(a,x)\big),\\
\big( b=b < (a,x)\big).
\end{cases}
\]
The first two elements in the second line lie in $A_{<a}$;  all other elements lie in $Z^\opp$.
The definition of the partial order on $M(f,\leq k-1)$ is such that we see that the maps define
together one from the product poset $(0<1)\times((A_{<a})^\opp\Join X_a)^\opp$ to 
$M(f,\leq k-1)$ and so are homotopic.
\\

\emph{Step 3: Conclusion.}
Let $\tilde Z\supset Z$  be the union of the mapping cones of the poset maps 
$\tilde s_a: (A_{<a})^\opp\Join X_a \to Z$ (these mapping cones have $Z$ in common).
Since $M(f,\leq k)$ is obtained from $M(f,\leq k-1)$ by putting a cone over each $A_{<a} *(X_a)$ with $\hgt (a)=k$, 
repeated use of the  above diagram yields a  homotopy equivalence 
$|\tilde Z| %{\buildrel\sim\over\longrightarrow}
\tilde \longrightarrow|M(f)|$.  Since $|A|$ is a deformation retract of  $|M(f)|$, 
we find a homotopy equivalence $|\tilde Z| %{\buildrel\sim\over\longrightarrow}
\tilde \longrightarrow|A|$.  Choose a 
homotopy inverse $H:|A| %{\buildrel\sim\over\longrightarrow}
\tilde \longrightarrow|\tilde Z|$.

Each composite  $\hat s_a=g\tilde s_a$ is null-homotopic by assumption (iv). This implies that the map $|g|$ 
extends to a continuous map $ G: |\tilde Z|\to |X|$. As $|\tilde Z|$ is  obtained from $|Z|$ by attaching 
cells of dimension $\ge n+1$, this map is still $n$-connected. 
% Again this attaching is up to homotopy, and an alternative argument looks at mapping cones.
So $GH$ is as desired.
\end{proof}

\section{The complex of unimodular decompositions}

We suppose $L$ \emph{unimodular} of genus $g>0$. We call a  subset $\ubold\subset\Ucal (L)_{>0}$ a 
\emph{unimodular decomposition} of $L$ if 
the natural map $\oplus_{u\in \ubold} u\to L$ is an isomorphism. If $\ubold$ and $\ubold'$ are 
unimodular decompositions, then we say that $\ubold'$ refines $\ubold$ (and we write $\ubold'\ge \ubold$) 
if every member of $\ubold'$ is contained in one of $\ubold$. This makes the collection of such 
decompositions a poset that we shall denote by $\Dcal (L)$. We chose this convention for the partial order  
(rather than its opposite) as to have  the standard height function on $\Dcal (L)$ assign to $\ubold$ 
the value $|\ubold|-1$ (the value $0$ being taken by the unique minimal element $(L)$). 
Notice that $\Dcal (L)$ is a subcategory of the category whose objects are finite subsets of 
$\Ucal (L)_{>0}$ and whose morphisms are surjections between these subsets. We sometimes write
$\Dcal_+ (L)$ for $\Dcal (L)_{>(L)}$, the subposet of \emph{strict} decompositions.

Recall that the \emph{barycentric subdivision} of a poset $X$ is the poset $X'$ whose elements are
chains in $X$ and for which $<$ is the relation `subchain of'. Its geometric realization is homeomorphic to that of $X$.

Now observe that there is a natural poset map
\[
f:(\Ucal (L)_{>0})'\to \Dcal (L)
\]
which assigns to a chain $0\not= u_0\subsetneq u_1\cdots\subsetneq u_k\subset L$ the decomposition 
whose members are $u_0, u_1\cap u_0^\perp, \dots ,u_k\cap u_{k-1}^\perp$ and
(in the case when $u_k\not= L$) $u_k^\perp $. 

\begin{theorem}\label{thm:dec}
The poset $\Dcal (L)$ of unimodular decompositions of $L$ is Cohen-Macaulay of dimension $g-1$.
\end{theorem}

We first prove a `simplicial counterpart' of this theorem.

\begin{lemma}\label{lemma:setdec}
Let $X$ be a finite set with at least two elements. Denote by $\Dcal_+(X)$ the poset of strict 
decompositions of $X$ (this excludes the trivial decomposition $\{ X\}$). Then $\Dcal_+(X)$ is spherical of dimension $|X|-2$.
\end{lemma}
\begin{proof}
We prove this with induction on $d:=|X|-2$. The case $d=0$ is trivial, so suppose $d>0$.
The poset $\Fcal (X)$ of proper nonempty subsets of $X$ is the 
the barycentric subdivision of the boundary of the simplex spanned by $X$ and is hence a combinatorial $d$-sphere. 
The poset $\Fcal (X)'$ of nested proper nonempty subsets of $X$ can be identified with the barycentric 
subdivision of $\Fcal (X)$ and so is still a combinatorial $d$-sphere.  
Now consider the  map $g_+: \Fcal (X)'\to \Dcal_+(X)$ which assigns to 
$\emptyset\subsetneq X_0\subsetneq X_1\subsetneq \cdots \subsetneq X_h\subsetneq X$ the partition 
$\{ X_0,X_1-X_0,\dots ,X_{h}-X_{h-1}, X-X_h\}$. This is a map of posets to which we want to apply 
Corollary \ref{cor:conn}. Let $\Ybold=\{ Y_0,\dots ,Y_{h+1}\}$ be a partition of $X$. Then
\[
\Dcal_+(X)_{> \Ybold}\cong \Dcal_+(Y_0)* \Dcal_+(Y_1)*\cdots*  \Dcal_+(Y_{h+1})
\]
and so by our induction hypothesis, $\Dcal_+(X)_{>\Ybold}$ is spherical of dimension  $-1+\sum_{i=0}^{h+1} \big(|Y_i|-1\big) =
|X|-3-h=d-1-h$. On the other hand, $g_+/\ubold$ can be identified with $\Fcal(\{ 0,\dots ,h+1\})'$, 
hence is a sphere  of dimension $h$. Now apply  version ($C^\opp$) of Corollary  \ref{cor:conn}.
\end{proof}

\begin{proof}[%
Proof 
of Theorem \ref{thm:dec}]
We proceed with induction on $g\ge 1$. For $g= 1$,  $\Dcal (L)=\{ L\}$ and so there is nothing to show. 
Assume $g>1$ and the theorem proved for genera $<g$. We apply version ($C^\opp$) of 
Corollary  \ref{cor:conn} to $f$. Let $\ubold=\{ u_0,\dots, u_{h}\}\in \Dcal(L)$, so that $\hgt (\ubold)=h$. 
Then we readily observe that
\[
\Dcal(L)_{>\ubold}\cong\Dcal_+(u_0)*\cdots*  \Dcal_+(u_{h})
\]
and so $\Dcal (L)_{>\ubold}$ is spherical of dimension  $-1+\sum_{i=0}^{h} \big(g(u_i)-1\big) =g-2-h$. 
On the other hand, $f_+/\ubold$ can be identified with the poset of nested nonempty subsets  of $\{ 0,\dots ,h\}$. 
This is contractible as it has the entire set as its maximal element.
So the conditions  ($C^\opp$) of \ref{cor:conn} are satisfied. Combining this  with 
Theorem \ref{thm:mc} yields that $\Dcal(L)$  is spherical of dimension $g-1$. In passing we showed 
that $\Dcal (L)_{>\ubold}$ is spherical. Notice that $\Dcal (L)_{<\ubold}$ can be identified with the poset 
$\Dcal_+(\{ 0,\dots ,h\})$ of all
strict decompositions of the set  $\{ 0,\dots ,h\}$. This is spherical of dimension $h-1$ by 
Lemma \ref{lemma:setdec} above. Similarly on finds that if $\ubold<\ubold'$, then 
$\Dcal(L)_{(\ubold,\ubold')}$ is spherical of dimension $|\ubold '|-|\ubold|-2$.
\end{proof}

The poset $\Dcal(L)$ appears in the moduli space of principally polarized abelian varieties. 
Consider 
the complexification of $L$, $L_\CC:=\CC\otimes_\ZZ L$ and extend $\la\, ,\, \ra$ bilinearly to $L_\CC$.  
On $L_\CC$ we also have the Hermitian form $H$ defined by $H(v,w):= \sqrt{-1}\la v ,\overline{w} \ra$, 
whose  signature is $(g,g)$. The space $\Scal (L)$ of complex subspaces $F\subset L_\CC$ of dimension 
$g$ on which $\la\, ,\, \ra$ is zero and $H$ is positive definite can be understood as the  
moduli space of complex structures on the real torus $(\RR/\ZZ)\otimes_\ZZ L$ that are polarized by $H$. 
It is a symmetric space for the symplectic group of $L_\RR$ and a choice of symplectic basis yields an 
isomorphism with the Siegel upper half space of genus $g$. Notice that  a unimodular decomposition  
$\ubold$ of $L$ determines a proper embedding of  $\prod_{u\in \ubold}\Scal (u)$
in $\Scal (L)$ with image a totally geodesic analytic submanifold. We shall denote that submanifold by $\Scal (\ubold)$.

\begin{corollary}\label{cor:dec}
The locus $\Scal (L)_\dec$ in the genus $g$ Siegel space $\Scal (L)$ which parametrizes the 
principally polarized abelian varieties that are decomposable as polarized varieties is a closed  
analytic subvariety of $\Scal (L)$ which has the homotopy type  of a bouquet of $(g-2)$-spheres.
In fact, the (closed) covering of $\Scal (L)_\dec$ by its irreducible components  is a Leray covering 
(which means that every nonempty finite intersection is contractible) whose nerve is identified with
$\Dcal_+(L)^\opp$ via the correspondence $\ubold\mapsto\Scal (\ubold)$. 
\end{corollary}
\begin{proof}
It is clear that $\Scal (L)_\dec$   is the union of the $\Scal (\{u,u^\perp\})$, where $u$ runs over the 
unimodular sublattices of $L$ that are neither $0$ nor $L$. It is known that this union is locally finite. 
So $\Scal (L)_\dec$ is a closed analytic subvariety of $\Scal (L)$ and as there are no inclusion relations among the 
$\Scal (\{u,u^\perp\})$, these yield the distinct  irreducible components of  $\Scal (L)_\dec$.  
In view of Weil's nerve theorem it now suffices to prove the  last statement.  

It is a  classical result (see for instance  \cite[\S 6.9]{debarre}) that any principally 
polarized abelian variety $A$  of positive dimension has a unique decomposition into indecomposables. 
Precisely,
if $\{A_i\}_i$ is the collection of abelian subvarieties of positive  dimension of $A$  that receive from 
$A$ a principal polarization and are  minimal for this property, then the natural map $\prod_i A_i\to A$ is 
an isomorphism. It follows  that  $\ubold\mapsto \Scal (\ubold)$ defines an bijection
from  $\Dcal_+(L)$ to the collection of nonempty intersections of the irreducible components of $\Scal (L)_\dec$.  
If we give the latter the poset structure defined by inclusion, then this is in fact a poset
isomorphism from $\Dcal_+(L)^\opp$. Each $ \Scal (\ubold)\cong\prod_{u\in\ubold}\Scal (u)$ is 
contractible and so the  covering has the asserted properties.
\end{proof}

\begin{corollary}\label{cor:agdec}
The locus $\Acal_{g,\dec}\subset \Acal_g$ of decomposables in the moduli space of principal 
polarized abelian varieties of  genus $g$ is a quasiprojective subvariety.
If $\tilde\Acal_g\to \Acal_g$ is a finite covering defined by a torsion free subgroup of 
$\Sp (2g,\ZZ)$ of finite index, and $\tilde \Acal_{g,\dec}$ denotes its preimage in $\tilde\Acal_g$, then the pair 
$(\tilde\Acal_g, \tilde \Acal_{g,\dec})$ is $(g-2)$-connected. In particular, $H_k(\Acal_g,\Acal_{g,\dec};\QQ)=0$ 
for  $k\le g-2$.
\end{corollary}
\begin{proof}[Sketch of Proof]
The first statement is in fact well-known, so let us just outline its proof.
It follows from Corollary \ref{cor:dec} that $\Acal_{g,\dec}\subset \Acal_g$ is a closed analytic subvariety. 
The Baily-Borel theory shows that its closure (relative to the Hausdorff topology) in the 
Baily-Borel compactification  $\Acal_g^*$ of  $\Acal_g$ is projective, as it is the image of an 
analytic morphism from the union of 
$\Acal_k^*\times \Acal_{g-k}^*$, $k=1,\dots, \lfloor \half g\rfloor$ to $\Acal_g^*$. Since $\Acal_g^*$ is projective,
and has a projective boundary, it follows that $\Acal_{g,\dec}$ is quasiprojective.

As to the remaining statements, let $\Scal_g$ be the Siegel space attached to the standard 
symplectic lattice $\ZZ^{2g}$.  According to Corollary \ref{cor:dec}, $\Scal_{g,\dec}$  has the 
homotopy type of a bouquet of $(g-2)$-spheres. So 
we can construct a relative CW complex $(Z,\Scal_{g,\dec})$ obtained from $\Scal_{g,\dec}$ by 
attaching cells of dimension $\ge g-1$ in a $\Sp(2g,\ZZ)$-equivariant manner as to ensure that 
$Z$ is contractible and no nontrivial element of $\Sp (2g,\ZZ)$  fixes a cell. 
If $\tilde\Acal_g$ is defined by the subgroup 
$\G\subset \Sp (2g,\ZZ)$, then $\G$ acts freely on $Z$ (as it does on the contractible $\Scal_g$) 
and so there is a $\G$-equivariant homotopy equivalence  $Z\to \Scal_g$ relative to  $\Scal_{g,\dec}$.  
It follows that there is also a homotopy 
equivalence $\G\bs Z\to \tilde\Acal_g$ relative to $\tilde \Acal_{g,\dec}$. 
Hence  $(\tilde\Acal_g, \tilde \Acal_{g,\dec})$ is $(g-2)$-connected.
\end{proof}

\begin{remark}
This corollary has a counterpart for the moduli space of stable genus $g$ curves with compact jacobian
(see \cite[Cor.\ 1.2]{looij}). 
\end{remark}

\section{Improved homological stability for the symplectic groups}

As an aside we show in this section that our main result yields a slight improvement of a result of Charney. 
We denote  the basis elements of $\ZZ^{2g}$ by $(e_1,e_{-1},\dots .e_g, e_{-g})$ and endow 
$\ZZ^{2g}$ with the symplectic form that this notation suggests: $\la e_i,e_{-i}\ra=\sign (i)$ and  
$\la e_i,e_j\ra=0$ when $j\not=\pm i$.  Let us write $G_g$ for the algebraic group attached to 
$\Sp(2g, \ZZ)$ so that for any ring $R$, 
$G_g(R)=\Sp (2g,R)$.  The obvious embedding of $G_g(R)$ in $G_{g+1}(R)$ (which identifies the 
$G_g(R)$ as the $G_{g+1}(R)$-stabilizer of the last two basis vectors)  induces a map on homology
that  is known to be an isomorphism in low degree for many choices of $R$:  
Charney \cite{charneyV}, following up on earlier work of Vogtmann, proved such a result for noetherian rings 
$R$ of finite noetherian dimension.  Her stability range is phrased in terms of the connectivity 
properties of a poset  $\HU_g(R)$  whose elements are what she calls \emph{split unimodular sequences} in 
$R^{2g}$: these are sequences of pairs $\big((v_1,v_{-1}),\dots , (v_k,v_{-k})\big)$ in $R^{2g}$  such that  
$\la v_i,v_{-i}\ra=\sign (i)$ and  $\la v_i,v_j\ra=0$ when $j\not=\pm i$
(so this is equivalent to giving a $k$ and a symplectic embedding of $\ZZ^{2k}$ in $\ZZ^{2g}$). 
The result may the be stated
as follows: if $a\ge \dim (R)+2$ is an integer such that $\HU_g(R)$ is $\lfloor\half(g-a-1)\rfloor$-connected 
for all $g$, then $H_i(G_g(R))\to H_i(G_{g+1}(R))$ is an isomorphism for $g\ge 2i+a+1$ and surjective for $g=2i+a$. 
This is supplemented by the
theorem that the connectivity assumption is fulfilled for the case $R$ is Dedekind domain resp.\  a 
principal ideal domain for $a=5$ resp.\ $a=4$. The point we wish to make is that we can do slightly 
better when $R$ is a Euclidean ring:

\begin{theorem}\label{thm:improvedstab}
If $R$ is a Euclidean ring, then we may take $a=2$:  $\HU_g(R)$ is $\lfloor \half(g-3)\rfloor$-connected and 
hence $H_i(G_g(R))\to H_i(G_{g+1}(R))$ is an isomorphism for $g\ge 2i+3$ and surjective for $g=2i+2$. 
\end{theorem}

There is corresponding statement with twisted coefficients (see Theorem 4.3 of  \emph{op.\ cit.}) 
that we will not bother to explicate.
The proof of Theorem   \ref{thm:improvedstab} rests on the observation that if we view the 
contents of the Theorems \ref{thm:mc} and \ref{thm:dec} as properties regarding $G_g(R)$ for $R=\ZZ$, 
then both statements and proofs remain valid if we let $R$ be any a Euclidean ring 
(that property enters in the proof of 
Theorem \ref{thm:maazen}, a result that is due to Maazen). And we may replace 
$\dim (R)$ with zero because  Charney uses $\dim (R)$ only to deal with projective modules. 
In our case they are free.

For the proof of \ref{thm:improvedstab} we need

\begin{lemma}\label{lemma:parts}
Let $A$ be a set endowed with a partition $P$ into a finite number of (nonempty) subsets. Then the poset 
 of sequences in $A$ which hit every part of $P$ at most once is spherical of dimension $|P|-1$.
\end{lemma}
\begin{proof}
Denote this poset by $X$ and let $Y$ be the allied poset of nonempty finite subsets of $A$ 
that meet every part at most once.
So we have an evident poset map  $f:X\to Y$. Since $|Y|$ can be identified with the iterated join of 
the distinct parts of $A$ (which are $|P|$ in number), it is spherical of dimension $|P|-1$.  
Observe that for any $y\in Y$, $f/y$ is the poset of nonempty sequences whose terms are distinct and lie in $y$. 
This is well-known to be spherical of dimension $|y|-1$ (\cite[Thm.\ 2.1]{maazen},  or 
\cite[Lemma 2.13(ii)]{vdk}). On the other hand, $|Y_{>y}|$ can be identified with the iterated join of the 
distinct parts of $A$ not hit by $y$ and as there are $|P|-|y|$ such parts, $|Y_{>y}|$ is
spherical of dimension $|P|-|y|-1$. Now apply Corollary \ref{cor:conn}.
\end{proof}

\begin{proof}[%
Proof 
of Theorem \ref{thm:improvedstab}]
Consider the poset map $f: \HU_g(R)\to \Dcal_+(R^{2g})$ which assigns to 
$\big((v_1,v_{-1}),\dots , (v_k,v_{-k})\big)$ the collection of genus 1 summands 
$U_i:=Rv_i+Rv_{-i}$, $i=1,\dots, k$ and the 
perp of their sum in $R^{2g}$. Given $\ubold\in \Dcal_+(R^{2g})$, denote by $t(\ubold )$ the number of genus 1 summands in 
$\ubold$ and  by $s(\ubold)$ the number of summands of higher genus.
Then clearly  $2s(\ubold)+t(\ubold)\le g$, or equivalently, $s(\ubold)\le \half(g-t(\ubold))$.
Observe that $f/\ubold$ is a poset of the type that appears in Lemma \ref{lemma:parts} above:
a member of $A$  is an oriented basis of any genus 1 summand in $\ubold$ and $P$ is the obvious partition defined by 
those summands. 
As there are $t(\ubold)$ parts, the lemma tells us that 
$f/\ubold$ is $(t(\ubold )-2)$-connected. On the other hand, $\Dcal_+(R^{2g})_{>\ubold}$
is a join of all the $\Dcal_+(u)$, $u\in \ubold$ and hence is connected of dimension 
$-2+\sum_{u\in\ubold} (g(u)-1)=g-2-|\ubold |$. In view of the inequality 
\[
t(\ubold )+ (g-|\ubold |)-2=g-s(\ubold)-2\ge \lceil\half(g+t(\ubold))\rceil-2\ge  \lfloor\half (g -3)\rfloor,
\]
it follows from the ($C^\opp$)-variant of Corollary \ref{cor:conn} that $\HU_g(R)$ is 
$\lfloor\half (g -3)\rfloor$-connected.
\end{proof}

\section{Relation to the separated curve complex}

A variation of the complex $\Dcal_+(L)$ considered above appears in the study of the Torelli groups.
Let $S$ be a closed connected orientable surface of genus $g$. 
Then $H_1(S)$ has the structure of a unimodular lattice of genus $g$ so that we have defined $\Dcal (H_1(S))$.
An isotopy class of embedded (unoriented) circles in $S$ is called \emph{separating}
if the complement of a representative decomposes $S$ into two connected components of positive genus.  
Notice that then the  homology of the components determines a nontrivial unimodular splitting of  $H_1(S)$. 
The \emph{separated curve complex} $\Ccal_\sep (S)$ of $S$ is the simplicial complex whose vertex set  are 
the isotopy classes of separating curves; a finite set of these spans a simplex precisely when its elements 
can be represented by embedded circles that are pairwise disjoint.  
The mapping class group $\G (S)$ ($=$ the group of connected components of the orientation preserving 
diffeomorphisms of $S$) acts on both $H_1(S)$ and $\Ccal_\sep (S)$.  The action on 
$H_1(S)$ has  the full integral symplectic group $\Sp(H_1(S))$ as its image; the kernel of this action is called the 
\emph{Torelli group} of $S$, denoted here by $T(S)$. The map that assigns to a separating curve its 
associated unimodular splitting of $H_1(S)$ extends
to a $\Sp(H_1(S))$-equivariant poset map  $T(S)\bs \Ccal_\sep (S)\to \Dcal_+(H_1(S))$. 
This is however not an isomorphism: think of the case when $g$ embedded circles split off 
$g$ genus 1 surfaces with a hole, so that what remains is a sphere with $g$ holes. Nevertheless:

\begin{theorem}\label{thm:curvecx}
The poset map $T(S)\bs \Ccal_\sep (S)\to \Dcal_+(H_1(S))$ is a homotopy equivalence, in particular, 
$T(S)\bs \Ccal_\sep (S)$ is $(g-3)$-connected.
\end{theorem}

This may be regarded as a natural companion of a result of one of us \cite{looij} that states $ \Ccal_\sep (S)$ is 
$(g-3)$-connected as well.

Before getting into the proof, we give a construction of $T(S)\bs \Ccal_\sep (S)$ entirely in terms of
$\Dcal_+(H_1(S))$.
\\

In what follows a \emph{tree}  is a finite simplicial complex of dimension $1$ that is contractible.
This amounts to giving a pair $T=(T_0,T_1)$, where $T_0$ is a finite set 
(the vertex set) and $T_1$ is a  \emph{nonempty} collection of two-element subsets of $T_0$ 
(the set of edges) such that  $|T_0|=|T_1|+1$ and its geometric realization is connected. 
The \emph{degree of a vertex}  $v\in T_0$, $\deg_T(v)$, is the number of elements of $T_1$ that contain $v$. 

We denote by $T_{00}\subset T_0$ the collection of vertices of degree $\le 2$. 
The identity $|T_0|=|T_1|+1$ is equivalent to $\sum_{v\in T_0} (deg_T(v)-2)=-2$.  
If $d_i(T)$ denotes the number of vertices of degree $i$, then this amounts to $\sum _{i\ge 3} (i -2)d_i= d_1-2$, 
which shows that if we fix a bound on $|T_{00}|=d_1+d_2$, then we have only finitely many isomorphism classes.

Notice that if $E$ is a  nonempty set of edges of $T_1$, and if we contract every connected component  of $T-E$, 
then we have formed a quotient tree $\pi^E:T\to T^E$ with the property that $\pi^E$ maps $E$ maps 
bijectively onto the edge set of $T^E$. 

\begin{lemma}\label{lemma:uniquecontract}
Let $E$ and $E'$ be nonempty sets of edges of a tree $T$ such that there exists an isomorphism
$h:T^E\cong T^{E'}$  with the property that $\pi^{E'}$ and $h\pi^E$ have the same restriction to $T_{00}$. Then $E=E'$.
\end{lemma}
\begin{proof}
We prove this with induction on $|T_1|$. For $|T_1|=1$ there is nothing to show and so suppose $|T_1|>1$.  
We prove that there is a terminal edge $e$ of $T$  on which $\pi^E$ and $\pi^{E'}$ coincide, i.e., for which either 
$e\notin E\cup E'$ or $e\in E\cap E'$. This suffices:  in the first case we can apply the induction 
hypothesis to the tree $T/e$, obtained by contracting $e$, and the images of $E$ and $E'$ in $T/e$. 
In the second case, $\pi^E(e)$ is a terminal edge with terminal vertex $\pi^E(v)$ and similarly for 
$\pi^{E'}(e)$ and $\pi^{E'}(v)$.
Since the isomorphism $h$ takes $\pi^E(v)$ to $\pi^{E'}(v)$, it will take $\pi^E(e)$ to $\pi^{E'}(e)$.
So here too the induction hypothesis can be invoked to the tree $T/e$ and the images of $E$ and $E'$ in $T/e$.

Choose a terminal vertex $v\in T_{00}$ and let $e=\{ v,v'\}\in T_1$ be the corresponding terminal edge. 
If $\deg (v')=2$ (so that $v'\in T_{00}$), then $\pi^E$ and $\pi^{E'}$ coincide on $e$.

So it remains to consider the case when $\deg (v')\ge 3$. This implies that $T-e$ is disconnected. 
Assume that $e\in E$. Then $e$ subsists as an end edge in $T^{E}$ and $\pi^E(v)$ is terminal in $T^{E}$. 
But then $\pi^{E'}(v)$ is also terminal in $T^{E'}$. If $e\notin E'$ (so that $v$ is contracted by $\pi^{E'}$), 
then this can only happen if all but one of the components of $T-e$ get contracted by $\pi^{E'}$. 
So if  $v''\in T_{00}$ is a terminal vertex of such a component, then $\pi^{E'}$ maps $v$ and $v''$ to the same vertex.
%the geodesic string in $T$ from $v$ to $v''$ is collapsed by $\pi^{E'}$ to a single edge. 
Then $\pi^E$ will have the same property. But $e$ lies on the geodesic string in $T$ from $v$ to $v''$ and so 
$\pi^E(v)\neq \pi^E(v')$, contradiction. So $e$ lies in $E'$ also.
%$E$ meets this string in $e$ only. We thus have found an edge (namely the last member of this string) that  is 
%neither in $E$ nor in $E'$.
\end{proof}

Let $\ubold$ be a finite set with at least two elements. We define a  \emph{$\ubold$-tree} 
as a tree $T$ endowed with a map  $i:\ubold\to T_0$ from to $\ubold$ to the vertex set $T_0$ of $T$ with the 
property that its image contains $T_{00}$. In view of the preceding observation, there are only finitely many 
isomorphism classes of $\ubold$-trees. In fact, the number of vertices not in $\ubold$ such a tree can have is at most
$|\ubold|-2$ and so its number of edges will be at most $2|\ubold|-3$. This bound is attained when
$i$ is injective, all vertices in the image of $i$ have degree 1 and all other vertices have degree 3.
Since a tree automorphism is  determined by its restriction to the set of terminal vertices, 
an $\ubold$-tree has no automorphism other than the identity.
We denote the set of isomorphism classes of $\ubold$-trees by $\Tcal (\ubold)$. We say that a $\ubold$-tree is 
\emph{strict} if the map $i:\ubold\to T_0$ is injective and we denote the corresponding subset of 
$\Tcal (\ubold)$ by $\dot\Tcal (\ubold)$.

We turn $\Tcal (\ubold)$ into a partially ordered set by stipulating that $(T',i')\le (T,i)$  if 
$T'$ is obtained from $T$ by means of a contraction map whose composite with $i$ yields $i'$. 
Lemma \ref{lemma:uniquecontract} shows that this contraction is already determined by the associated injection 
$T'_1\subset T_1$. So
$\Tcal (\ubold)_{\le (T,i)}$ can be identified with the poset of nonempty subsets of $T_1$ and hence
$|\Tcal (\ubold)_{\le (T,i)}|$ with the simplex $\sigma (T)$.  By the above computation,  the dimension of 
$\Tcal (\ubold)$ is at most $2|\ubold |-4$. 

\begin{lemma}
The poset   $\Tcal (\ubold)$ is contractible.
\end{lemma}
\begin{proof}
If $\ubold$ has just two elements, then $\Tcal (\ubold)$ is a singleton (it has a single edge) and
that is certainly contractible. Suppose therefore that $|\ubold |\ge 3$. Then the cone 
$C(\ubold)$ on $\ubold$ (which has one vertex $*$ not in $\ubold$ and
edges $\{ *,u\}$, $u\in\ubold$) is clearly a $\ubold$-tree.
Notice that $\Tcal (\ubold)_{\le  C(\ubold)}$ is the simplex on 
$\ubold$ and hence is contractible. We also have a deformation retraction of 
$\Tcal (\ubold)\to\Tcal (\ubold)_{\le  C(\ubold)}$ which assigns to  every $\ubold$-tree
in $\Tcal (\ubold)$ the $\ubold$-tree in $\Tcal (\ubold)_{\le  C(\ubold)}$ obtained  by contracting all the
internal edges, that is, the edges that do not have a vertex of degree one. 
\end{proof}

We return to the poset $\Dcal_+(L)$ of strict unimodular decompositions of $L$. Observe that any 
$\ubold\in\Dcal_+(L)$ is a  finite subset of the set of all unimodular lattices $u$ with $0\not= u\not= L$ and that 
the relation $\ubold'\le \ubold$ determines (and is given by) a surjection $\ubold\to \ubold'$. 
So  the poset $\Dcal_+(L)$ is a subcategory of the category of sets.
We define a poset $\Tcal\Dcal (L)$ of \emph{tree decompositions} of $L$: 
its underlying set is the disjoint union of the  isomorphism classes of strict $\ubold$-trees, where $\ubold$ runs 
over $\Dcal_+(L)$:
\[
\Tcal\Dcal  (L):=\bigsqcup_{\ubold \in\Dcal (L)} \dot\Tcal (L)
\]
and we stipulate that $(T', i': \ubold'\hookrightarrow  T'_0)\le (T, i: \ubold\hookrightarrow  T_0)$ if 
$\ubold$ refines $\ubold'$  and there exists an edge contraction $T\to T'$ extending the surjection
 $\ubold\to\ubold'$ associated to the refinement property. It follows from Lemma \ref{lemma:uniquecontract} that this edge 
contraction is unique, so that we have defined a partial order indeed.
Notice that there is an obvious forgetful map of posets
\[
p: \Tcal\Dcal  (L)\to \Dcal_+(L).
\]

\begin{proposition}\label{prop:treeeq}
The poset map $p:\Tcal\Dcal(L)\to \Dcal_+(L)$ is a homotopy equivalence.
\end{proposition}
\begin{proof}
For every $\ubold\in \Dcal_+(L)$, $p/\ubold $ may be identified with $\Tcal (\ubold)$, hence is contractible. 
Now apply Proposition (1.6) of \cite{quillen}.
%Notice that $\Dcal_+(L)_{>\ubold}$ is spherical of dimension $g-3-\hgt (\ubold)$ and that we could also apply 
%Corollary \ref{cor:conn}.
\end{proof}

\begin{proof}[%
Proof 
of Theorem \ref{thm:curvecx}]
A simplex of $\Ccal_\sep (S)$ is given by a closed one-dimensional submanifold $A\subset S$ such that every 
connected component of 
$S-A$ has negative Euler characteristic and the associated graph $T(S,A)$ (having the set of connected components of 
$S-A$ resp.\ $A$ as its set of vertices resp.\ edges) is a tree. Since a connected component of $S-A$ of zero genus has 
at least 
$\ge 3$ boundary components, all vertices of degree $\le 2$ correspond with connected components of $S-A$ with
nonzero genus, hence positive genus. 
Moreover, the first homology groups of the connected components of positive genus define a decomposition of $H_1(S)$.
This gives $T(S,A)$ the structure of a tree decomposition of $H_1(S)$. Another choice $A'$ for $A$ yields the 
same tree decomposition of $H_1(S)$ if and only if there exists an orientation preserving self-homeomorphism of 
$S$ that takes $A$ to $A'$ and  induces the identity map in $H_1(S)$. It is now easy to identify the poset map
$p$ in Proposition 
\ref{prop:treeeq} with the poset map $T(S)\bs \Ccal_\sep (S)\to \Dcal_+(H_1(S))$ introduced
at the beginning of this section. The theorem follows.
\end{proof}

\section{Appendix}
In this appendix we recall the proof of the theorem of Maazen on the 
Cohen-Macaulayness of the poset 
$\Ocal(n)$ of partial bases in a free module $R^{n}$ over a Euclidean ring.

Given a partial basis $w=(w_{1},\dots, w_{r})$, we write $\Ocal(n)_{w}$
for the poset of partial bases $(v_{0},\dots,v_{s})$ so that 
$(v_{0},\dots,v_{s},w_{1},\dots, w_{r})$ is a partial basis of $R^{n}$.
Observe that $\GL(n)$ acts transitively on partial bases of a given length.
We write $p_{n}:R^n\to R$ for taking the last coordinate. For $a, b\in R$ with $||b||>0$, let $q(a,b)\in R$ be the scalar 
given by the Euclidean structure, so that $||a-q(a,b)b||<||b||$.
For $k\geq0$ the subposet $\Ocal(n,k)$ of $\Ocal(n)$
consists of the $(v_{0},\dots,v_{s})$ with $||p_{n}(v_{i})||\leq k$ for all $i$.
So $\Ocal(n,0)\cong \Ocal(n-1)$.

If $w=(w_{1},\dots, w_{r})\in  \Ocal(n)$ and $||p_n(w_{i})||=k$ for some $i\leq r$,
we define  $\rho_{w,i}:R^{n}\to R^{n}$ by $\rho_{w,i}(v)=v-q(p_{n}(v),p_n({w_{i}}))w_{i}$. This induces a 
retract
from 
$\Ocal(n)_{w}$ to $\Ocal(n,k-1)_w:=\Ocal(n,k-1)\cap \Ocal(n)_{w}$.
The key observation is that if $\Ocal(n)_{w}$ is $d$-connected, then so is
$\Ocal(n,k-1)_{w}$. The other tools that Maazen uses
have been 
published in \cite{charney}  and \cite{vdk}, to which we refer freely.

\begin{theorem}[Maazen]\label{thm:maazen} 
Let $R$ be a Euclidean ring and $d\in \ZZ$.
\begin{enumerate}\setlength{\itemindent}{1em}
\item[(i)]  $\Ocal(n)$ is $d$-connected if $n\geq d+2$,
\item[(ii)] $\Ocal(n)_{w}$ is $d$-connected for $w\in \Ocal(n)$
if $n\geq d+|w|+2$,
\item[(iii)] $\Ocal(n,k)_{w}$ is $d$-connected for $w\in \Ocal(n)$, if there is an $i$ with $||p_{n}(w_{i})||=k+1\geq1$
and $n\geq d+|w|+2$,
\item[(iv)] $\Ocal(n,k)$ is $d$-connected 
if $n\geq d+2$ and $k\geq1$.
\end{enumerate}
\end{theorem}

\begin{proof}
By induction on $d$. For $d<-1$ there is nothing to prove.
Let $n\geq d+2$. Choose $w\in \Ocal(n,1)\bs \Ocal(n,0)$ with $|w|=1$.
Then $p_{n}(w_{1})$ is a unit, so that the link of $w$ in
$\Ocal(n,0)$ is all of $\Ocal(n,0)$.
Apply
\cite[Lemma 2.13(ii)]{vdk} to conclude that $\Ocal(n,1)$ is $d$-connected. 
By induction on $k$ it follows from \cite[Lemma 2.13(i)]{vdk} that $\Ocal(n,k)$ is $d$-connected for 
$k\geq1$. Taking the limit $k\to\infty$ we see that $\Ocal(n)$ is $d$-connected. 
For case (ii) we may assume that $w$ consists of the last $|w|$ elements of the standard basis.
By \cite[Cor.\ 1.7]{charney}, with $S$ equal to the span of $w$ and $F=\Ocal(n-|w|)$,  case (ii) follows.
Then case (iii) follows by means of the retract.
\end{proof}

{}From \cite[Cor.\ 1.7]{charney} with $n=|w|$ we see that $\Ocal(m)_{>w}$ is $(m-|w|-2)$-connected for $w\in \Ocal(m)$.
One concludes as in \ref{prop:st} that

\begin{corollary}
$\Ocal(n)$ is Cohen-Macaulay.
\end{corollary}

\end{document}